\documentclass[a4paper,notitlepage]{siamsials251208}

\usepackage{amsfonts}
\usepackage{graphicx}
\usepackage{epstopdf}
\usepackage{algorithmic}
\ifpdf
  \DeclareGraphicsExtensions{.eps,.pdf,.png,.jpg}
\else
  \DeclareGraphicsExtensions{.eps}
\fi

\usepackage{enumitem}
\setlist[enumerate]{leftmargin=.5in}
\setlist[itemize]{leftmargin=.5in}

\newsiamremark{remark}{Remark}
\newsiamremark{hypothesis}{Hypothesis}
\crefname{hypothesis}{Hypothesis}{Hypotheses}
\newsiamthm{claim}{Claim}
\newsiamremark{fact}{Fact}
\crefname{fact}{Fact}{Facts}

\headers{Emerging Thermal-Metabolic Control}{M. A. Neto, P. A. Marquet, M. A. Freilich, L. Martí, and N. Sanchez-Pi}

\title{Emerging Thermal-Metabolic Control of Phytoplankton Annual Cycles: Extending the Critical Depth Hypothesis%
\thanks{Submitted to SIAM Journal on Life Sciences
\funding{This work was supported by the Franco-Chilean Binational Center of Artificial Intelligence, ANID Strengthening R\&D capabilities Program CTI230007 Inria Chile, and Inria Project MOFOMO - Multimodal Ocean Foundation Models for a Digital Twin of the Ocean. P.A.M. acknowledges support from Centro de Modelamiento Matemático (CMM), Grant FB210005, BASAL funds for Centers of excellence from ANID-Chile, the International Center for
Theoretical Physics (ICTP) through their Associates Programme funded by the Simons Foundation through grant number 284558FY19, and the Earth Comission and Global Commos Alliance.}}}

\author{Matías A. Neto\thanks{Inria Chile Research Center, Las Condes, RM, Chile 
  (\email{matias.neto@inria.cl}, \email{luis.marti@inria.cl}, \email{nayat.sanchez-pi@inria.cl}).}
\and Pablo A. Marquet\thanks{Facultad de Ciencias Biol\'ogicas, Pontificia Universidad Católica de Chile, Santiago; The Santa Fe Institute, NM, USA; Centro de Modelamiento Matemático (CMM), Chile; Centro de Cambio Global UC, Chile; Instituto de Sistemas Complejos de Valpara\'iso (ISCV), Chile; and  Wissenschaftskollege zu Berlin, Berlin, Germany (\email{pmarquet@uc.cl}).}
\and Mara A. Freilich\thanks{Division of Applied Math and Department of Earth, Environmental, and Planetary Sciences, Brown University, Providence, RI, USA (\email{mara\_freilich@brown.edu}).}
\and Luis Martí\footnotemark[2]
\and Nayat Sanchez-Pi\footnotemark[2]}

\usepackage{amsopn}

\ifpdf
\hypersetup{
  pdftitle={The Metabolic and Physical Drivers of Phytoplankton Persistence in a Changing Ocean},
  pdfauthor={M. A. Neto, P. A. Marquet, M. A. Freilich, L. Martí, N. Sanchez-Pi}
}
\fi

\usepackage{xspace}

\usepackage[version=4]{mhchem}
\usepackage{graphicx}%
\usepackage{subcaption}
\usepackage{multirow}%
\usepackage{xcolor}%
\usepackage{textcomp}%
\usepackage{manyfoot}%
\usepackage{booktabs}%
\usepackage{listings}%
\usepackage{rotating}
\usepackage{comment}
\usepackage{siunitx}
\usepackage{multiobjective}
\usepackage{xfrac}
\usepackage{makecell}
\usepackage{cite}
\usepackage[mathlines]{lineno}%
\usepackage[showdeletions]{color-edits}
\addauthor[Matias]{mn}{orange}
\addauthor[Pablo]{pm}{magenta}
\addauthor[Nayat]{nspi}{BrickRed}
\addauthor[Luis]{lm}{cyan}
\addauthor[Mara]{mf}{ForestGreen}

\usepackage[T1]{fontenc} 
\newcommand{\etal}{\emph{et\,al.}\xspace}
\newcommand{\GFDLHyperThreshold}{0.80}

\newcommand{\GFDLHyperAboveThresholdPerc}{83.1}
\newcommand{\GFDLHyperMeanRsq}{0.859}

\begin{document}
\maketitle

%%==================================%%
\begin{abstract} While the Critical Depth Hypothesis (CDH) established physical mixing and light availability as primary drivers of phytoplankton annual cycle, global thermal restructuring requires expanding this paradigm to incorporate temperature-dependent biological kinetics. We introduce a non-autonomous framework that extends the CDH into a temperature-dependent paradigm. By analyzing invasion dynamics at the boundary of extinction, we derive the critical nutrient saturation($\gamma_\text{crit}$) as a structural threshold for population persistence across annual cycles. Applied to end-of-the-century projections from the GFDL-ESM4 model under a high-emission scenario (SSP5-8.5), our framework identifies an emerging thermal-metabolic control, which increasingly displaces regions where seasonal dynamics were historically governed by physical mixing.
\end{abstract}
\begin{relevance}
Disentangling the primary drivers of phytoplankton annual cycles is essential for understanding ecosystems under climate change. Our analysis extends CDH by demonstrating that thermal metabolism modifies the baseline required for annual population recovery. Projections under high-emission scenarios reveal that anthropogenic forcing alters this balance, shifting global phytoplankton driver control from physical mixing to thermal-metabolic. However, this transition is not homogeneous: the North Atlantic Subpolar Gyre serves as a mixing-driven anomaly where deep mixing dynamics and regional thermal stasis maintain historical control. By evaluating stability directly from model outputs, this methodology provides a computationally efficient diagnostic to map driver hierarchies across Earth System Models.
\end{relevance}
\begin{mathcontent}
We model seasonal ecosystem dynamics using a non-autonomous NPZ-t framework forced by periodic functions of temperature, light, and mixed layer depth. Linearizing around the extinction equilibrium $E_0$ decouples the system into a diagonal periodic structure. The Floquet multiplier $\rho_P$ yields the annual invasion growth rate $\lambda_P$, we prove that $\lambda_P>0$ is a sharp condition for uniform persistence (Theorem \ref{th:persistence}). Isolating the threshold condition $\lambda_P=0$ yields the definition of the critical nutrient saturation ratio, $\gamma_\text{crit}$. A first-order Taylor sensitivity decomposition of $\gamma_\text{crit}$ allow us to assess the role of environmental anomalies for temperature ($\Delta T$), mixing ($\Delta H$), and light amplitude ($\Delta I$) on phytoplankton annual cycles.
\end{mathcontent}

% REQUIRED

% REQUIRED



% \pagestyle{fancy}
% \fancyhead{} % clear all header fields
% \renewcommand{\headrulewidth}{0pt} % no line in header area
% \fancyfoot{} % clear all footer fields
% \fancyfoot[LE,RO]{\thepage}           % page number in "outer" position of footer line
% \fancyfoot[RE,LO]{\textcolor{red}{\footnotesize Work in progress, do not disseminate. Comments welcomed!}}

\section{Introduction}\label{sec1}

Marine microbial communities play a central role in regulating atmospheric and oceanic chemistry through photosynthesis and carbon sequestration \cite{falkowski2008microbial,fuhrman2015marine,litchman2012phytoplankton}. Phytoplankton annual cycles have long been framed around two primary drivers: light limitation and mixed layer depth (MLD), as formalized by Sverdrup's Critical Depth Hypothesis (CDH) \cite{sverdrup}. However, the global ocean is undergoing sustained but heterogeneous thermal change, ranging from accelerated warming in polar regions \cite{meredith2019polar} to anomalous cooling in subpolar gyres \cite{rahmstorf_2015}. Temperature modulates biological rates across all scales of biological organization \cite{hochachka2002biochemical,brown2004,dell2014}, modifying growth and mortality rates in plankton populations \cite{Eppley1972,Robarts1987}. Extending the CDH light-mixing framework to incorporate thermal metabolism impacts is essential to understand phytoplankton annual dynamics in a changing ocean.

Current Earth System Models (ESMs) provide essential projections of ecosystem responses, yet their structural complexity often obscures the mechanistic hierarchy of drivers controlling annual biological cycles \cite{anderson_2005}. Disentangling drivers of change in marine ecosystems is important for improving confidence in predictions and gaining mechanistic understanding about changes in ecosystems that result from physical forcing \cite{Rohr2023,martin_when_2024}. Marine ecosystems are highly variable with timescales of variability ranging from diel to seasonal to decadal. The theory of non-autonomous periodic systems provides a mathematical foundation to study change in the context periodically (seasonally) forced dynamics \cite{zhao,wang2008threshold} and has been explored in specific ecological models \cite{MingChen2017,klausmeier_floquet_2008}. However, an analytical bridge capable of disentangling driver dominance across global ocean models remains absent.

Here, we extend the CDH into a temperature-dependent paradigm to diagnose the primary drivers of phytoplankton annual cycles in a changing climate. By analyzing invasion dynamics, we derive the critical nutrient saturation, $\gamma_{\text{crit}}$, metric that represent the baseline required for annual population recovery. Applying a sensitivity decomposition to global ESM projections (GFDL-ESM4 under SSP5-8.5), we define realized impacts of thermal, physical and light drivers. This framework reveals a global regime shift in annual cycle controls, where thermal metabolism displaces physical mixing as the primary driver of phytoplankton dynamics on seasonal timescales.

\section{NPZ-t: Non-autonomous NPZ model} \label{sec:threshold}

Let $N(t)$, $P(t)$, and $Z(t)$ denote, respectively, the vertically averaged nutrient, phytoplankton, and zooplankton concentrations within the mixed layer at time $t\geq0$ within a model grid cell (at a specific horizontal location). Assuming a closed inventory $C_0>0$ and a period $T_p>0$, we define the NPZ-t dynamics as
\begin{align}\label{eqn:NPZ-t-full}
&\left(N(0),P(0),Z(0)\right)\in X:= \left\{(N,P,Z)\in\mathbb{R}^3_+: N+P+Z = C_0 \right\},\\
&\frac{dN}{dt}=-f(t)\mu_P(t) g(N)P+(1-\alpha)\mu_Z(t)h(P)Z+m_P(t)P+m_Z(t) Z^2+s(t)(P+Z),\\
&\frac{dP}{dt}=f(t)\mu_P(t)g(N)P-\mu_Z(t)h(P)Z-m_P(t)P-s(t)P,\\
&\frac{dZ}{dt}=\alpha \mu_Z(t)h(P)Z-m_Z(t) Z^2-s(t)Z,
\end{align}
where $f(\cdot)$ represents light absorption, $g(\cdot)$ captures the nutrient uptake, $h(\cdot)$ expresses the zooplankton grazing. Moreover, $\mu_P$ and $\mu_Z$ denote maximum growth rate and grazing rate, respectively. Also, $m_P(\cdot)$ and $m_Z(\cdot)$ denote mortality rates, $s$ accounts for dilution from entrainment of water below the mixed layer, and $\alpha\in (0,1)$ denotes the zooplankton assimilation efficiency. This is a bulk mixed layer model, which models average dynamics within the mixed layer and neglects lateral gradients \cite{evans,Mara}.

We will assume that $f,s,\mu_P, \mu_Z, m_P,m_z$ are non-negative, $C^0$ and $T_p$-periodic functions. To ensure the periodic renewal of the mixed layer, we assume strictly positive annual mixing ($\int^{T_p}_0 s(t)dt> 0$). The growth and grazing functions $g$ and $h$ are non-negative and of class $C^1$. The functionals forms obey the assumptions
\begin{align*}
    \tag{H1}\label{hyp:H1}
   g(N) \le g(C_0), \quad \forall N \in [0, C_0], \\
   \tag{H2}\label{hyp:H2}
    \exists\lim_{P\to 0^+}\frac{h(P)}{P}=:h_0\in[0,\infty).
\end{align*}
Assumption \eqref{hyp:H1} enforces that nutrient uptake is bounded by its maximal capacity $g(C_0)$. Assumption \eqref{hyp:H2} guarantees that per-capita grazing pressure remain bounded at low phytoplankton densities and vanishes continuously near the origin.

\subsection{Invasion Analysis}

Following modern invasion theory \cite{Chesson2000}, we evaluate phytoplankton long-term persistence throughout the annual cycle by analyzing the local linear stability of the extinction state $E_0 = (C_0, 0, 0)$, which represents the unexploited resource state.

By mass conservation on the invariant simplex $X$, we substitute $N = C_0 - P - Z$ to reduce the dynamics to the $(P, Z)$ subspace. Linearizing around $(P, Z) = (0, 0)$ and noting that assumption \eqref{hyp:H2} implies $h(0) = 0$, the system decouples into a diagonal periodic system

\begin{equation}\label{local-system}
    \frac{d}{dt} \begin{pmatrix} P \\ Z \end{pmatrix} = J(t) \begin{pmatrix} P \\ Z \end{pmatrix}, \quad
    J(t) = 
    \begin{bmatrix}
      f(t)\mu(t)g(C_0) - m_P(t) - s(t) & 0 \\
      0 & - s(t)
    \end{bmatrix}\,.
\end{equation}
According to Floquet theory \cite{Hale}, the stability of $E_0$ depends on the eigenvalues of the monodromy matrix $C = \Phi(T_p)$. Direct integration of the decoupled system \eqref{local-system} yields the explicit matrix
\begin{equation*}
    C = 
    \begin{bmatrix}
    \exp\left(\int_0^{T_p} [f(\tau)\mu_P(\tau)g(C_0)-m_P(\tau)-s(\tau)] \, d\tau\right) & 0\\
    0 & \exp\left(-\int_0^{T_p} s(\tau) \, d\tau\right)
    \end{bmatrix}
\end{equation*}
where the diagonal elements correspond to the Floquet multipliers for phytoplankton, $\rho_P$, and zooplankton, $\rho_Z$, respectively. As we assume strictly positive annual mixing, the zooplankton multiplier always satisfies $|\rho_Z| < 1$.

Thus, the stability of the extinction state depends entirely on the phytoplankton multiplier $\rho_P$. This relationship defines the stability threshold via the annual invasion growth rate $\lambda_P$
\begin{equation}
  \lambda_P = \frac{1}{T_p} \ln(\rho_P) = \frac{1}{T_p}\int_0^{T_p} \left( f(\tau)\mu_P(\tau)g(C_0) - m_P(\tau) - s(\tau) \right) \, d\tau\,.  
\end{equation}
For a detailed mathematical foundations on the concepts underpinning Floquet theory and uniform persistence, we refer the reader to Appendix \ref{sec:foundations}.

\subsection{Global Extinction and Persistence Theorems}

We prove that the extinction state acts as a persistence threshold: If $E_0$ is locally stable, the ecosystem undergoes global extinction; when $E_0$ is locally unstable, a positive invasion growth rate repels the phytoplankton population from the zero boundary, ensuring uniform persistence. As formalized in the following theorems.

\begin{theorem}[Global Extinction Condition]\label{th:extinction}
Consider the NPZ-t system \eqref{eqn:NPZ-t-full}. If $E_0$ is locally stable, \emph{i.\,e.}, $\lambda_P<0$, then $E_0$ is globally asymptotically stable, \emph{i.\,e.},
\begin{equation}
  \lim_{t\to \infty} P(t)= 0 \text{ and } \lim_{t\to \infty} Z(t)= 0\,.  
\end{equation}

\end{theorem}

Theorem \ref{th:extinction} establishes that when $\lambda_P<0$, metabolic losses and physical dilution outpace the maximum potential growth even under saturating nutrient conditions, forcing the inevitable collapse of the phytoplankton population. Consequently, the absence of prey drives the secondary extinction of the zooplankton population.

\begin{theorem}[Global Persistence Condition]\label{th:persistence}
Consider the NPZ-t system \eqref{eqn:NPZ-t-full}. If $E_0$ is locally unstable, \emph{i.\,e.}, $\lambda_P>0$, then the system is uniformly persistent within the invariant set defined by $C_0$. Specifically, there $\exists~\eta>0$, such that all solutions satisfy
\begin{equation}
    \liminf_{t\to \infty}{P(t)}\geq \eta\,.
\end{equation}

\end{theorem}

The condition $\lambda_P>0$ ensures that the phytoplankton population repels from the zero boundary, guaranteeing long-term phytoplankton persistence throughout the annual cycle.

%\begin{remark}
%    The correspondence between the local instability of $E_0$ (whenever $\lambda_P>0$) an the global property of uniform persistence is consistent with the general threshold theory of periodic compartmental models \cite{wang2008threshold}. Specifically, \cite{MingChen2017} demonstrate that for seasonal systems, the condition $\rho(\Phi)>1$ -- where $\rho(\Phi)$ is the spectral radius of the evolution operator -- is equivalent to an ecological reproductive index $R_0>1$, which ensures the transition from extinction to a stable periodic bloom. In our framework, the structural decoupling at the extinction state -- where the nutrient remains at a constant baseline $C_0$ -- ensures that $\lambda_P$ is determined solely by the annual reproductive success against a stable resource background. This formally anchors the invasion growth rate $\lambda_P$ in the theory of periodic dynamical systems, providing a rigorous generalization of the classical Critical Depth criteria to non-autonomous nutrient-limited environments.
%\end{remark}

\subsection{The Critical Nutrient Saturation}

Theorems \ref{th:extinction} and \ref{th:persistence} provide a binary criterion for phytoplankton persistence throughout the annual cycle based on the sign of the invasion exponent $\lambda_P$. However, evaluating this exponent globally presents a practical challenge: while physical variables are relatively well-constrained, the specific nutrient inventory $C_0$ remains relatively uncertain across global climate projections \cite{kwiatkowski, seferian2020tracking}. 

To overcome this limitation, we adopt an inverse nutrient approach. Instead of assuming a fixed nutrient load, we isolate the boundary condition where the invasion growth rate vanishes ($\lambda_P=0$). This allows us to calculate the minimum annual nutrient uptake efficiency required to prevent extinction.

By setting $\lambda_P = 0$ and separating growth gains from metabolic and physical losses, the threshold balance is formulated as
\begin{equation}\label{eq:balance}
    \gamma_\text{crit} \int_0^{T_p} \mu_P(t)f(t) \, dt = \int_0^{T_p} \left( m_P(t) + s(t) \right) \, dt\,,
\end{equation}
where $\gamma_\text{crit}$ represents the non-dimensional critical nutrient saturation required to prevent extinction. Defining the integrated annual growth capacity as $\mathcal{G}(x)$ and the integrated annual losses as $\mathcal{L}(x)$, we solve explicitly for $\gamma_\text{crit}(x)$
\begin{equation}\label{eq:gamma_crit}
    \gamma_\text{crit}(x) = \frac{\mathcal{L}(x)}{\mathcal{G}(x)}\,.
\end{equation}
If $\gamma_\text{crit}(x)=0.5$, the phytoplankton population must operate at $50\%$ of its maximum uptake capacity just to survive, any further limitation would lead to extinction. Hence, phytoplankton annual cycles are only possible if the total inventory allows for a saturation level that exceeds the threshold $g(C_0)>\gamma_\text{crit}$.

\subsection{Trivariate Driver Decomposition of $\gamma_\text{crit}$}

To identify the primary driver of phytoplankton annual cycles, we decompose the variability of $\gamma_\text{crit}$ into its thermal-metabolic, physical-mixing and light-amplitude drivers through local sensitivity analysis. 

Let $y$ denote a specific year (e.g., $2024$ or $2100$) and $\mathcal{P}_{\text{ref}}$ the historical reference period (1960--2000). We define the environmental anomalies for temperature ($\Delta T$), mixing ($\Delta H$), and light amplitude ($\Delta I$) as
\begin{align}
    \Delta T(x) &= \bar{T}_{y}(x) - \langle \bar{T}(x)\rangle_{\text{ref}}\,,\\
    \Delta H(x) &= H^{\max}_{y}(x) - \langle H^{\max}(x)\rangle_{\text{ref}}\,,\\
    \Delta I(x) &= I_{\text{amp}, y}(x) - \langle I_{\text{amp}}(x)\rangle_{\text{ref}}\,,
\end{align}
where $\bar{T}_{y}$ is the annual mean temperature, $H^{\max}_{y}$ is the annual maximum mixed layer depth, and $I_{\text{amp}, y}$ is the annual seasonal irradiance amplitude. We scale these anomalies by the local sensitivity gradients of the annual cycle threshold to calculate the realized impacts ($\mathcal{I}$)
\begin{equation}
    \mathcal{I}_{T}(x) = \left| \frac{\partial \gamma_\text{crit}}{\partial \bar{T}}\right| |\Delta T(x)|\,, \quad 
    \mathcal{I}_{H}(x) = \left| \frac{\partial \gamma_\text{crit}}{\partial H^{\max}}\right| |\Delta H(x)|\,, \quad
    \mathcal{I}_{I}(x) = \left| \frac{\partial \gamma_\text{crit}}{\partial I_{\text{amp}}}\right| |\Delta I(x)|\,.
\end{equation}
As $\gamma_\text{crit}$ is a dimensionless metric, the product of its local gradients and the corresponding environmental anomalies yields dimensionless realized impacts. This dimensional consistency allows us to directly compare the relative magnitudes of thermal, physical, and light drivers, establishing the quantitative basis to determine the driver control of the phytoplankton annual cycle.

\section{Global Ocean Implementation and Experimental Setup}

\subsection{Physical Forcing and Seasonal Dynamics} 
Following the approach of Evans and Parslow \cite{evans} and Freilich~\etal \cite{Mara}, the average light limitation within the MLD and the dilution effects due to deepening are modeled as
\begin{align}
f(t) & = \frac{1}{K_d H(t)}\left(1 - e^{-K_d H(t)}\right) \frac{I(t)}{I(t)+I_0}\,,\\
s(t) & = 
\begin{cases}
\frac{1}{H(t)} \frac{d H(t)}{dt}, & \text{if }\frac{d H(t)}{dt} \ge 0,\\
0, & \text{otherwise.}
\end{cases},
\end{align}

\noindent where $I(t)$ is the surface irradiance, $H(t)$ is the MLD, $T(t)$ is the Sea Surface Temperature (SST), $K_d$ is the light attenuation coefficient, $I_0$ is the light half-saturation constant and $s(t)$ captures the dilution of biomass as the mixed layer deepens. 

%In open-system formulations, entrainment introduces nutrients from an infinite deep reservoir via $s(N_{\text{deep}}-N)$; the total concentration converges to the deep reservoir concentration $N_\text{deep}$ as $t\to \infty$. Since this study focuses on long-term uniform persistence, fixing a closed inventory $C_0$ via instantaneous recycling $s(P+Z)$ allows the model to bypass the transient dynamics. While exact vertical integration of photosynthesis (\emph{e.\,g.} Huisman and Weissing \cite{Huisman}) provides detailed depth-resolved profiles, this bulk-mixed layer formulation captures the fundamental seasonal light limitation required for global-scale ESM coupling.

We approximate $I(t)$, $H(t)$, and $T(t)$ using periodic functions that capture the minima and maxima of the seasonal cycle. These variables are derived by fitting the periodic functions to seasonal climatologies from GFDL-ESM4 model \cite{gfdl-esm4} at each grid point. To focus on ice-free pelagic dynamics, we exclude grid points where the minimum annual temperature falls below $\qty{-1.8}{\unit{\degreeCelsius}}$. In these regions, cryospheric processes and extreme light limitation -- not explicitly modeled here -- significantly alter the biological response.

\subsection{Nutrient Uptake, Grazing and Metabolic Rates}

Nutrient limitation is modeled using Monod kinetics, while zooplankton grazing follows a Holling Type III functional response, consistent with the dilution-recoupling mechanism \cite{Mara,behrenfield2010}, as
\begin{align}
g(N) & = \frac{N}{N + N_0}\,,\\
h(P) & = \frac{P^2}{P^2 + P_0^2}\,,
\end{align}
where this functional forms are multiplied by $\mu_P$ and $\mu_Z$ that denote the time-dependent maximum growth and grazing rates, respectively, and $N_0$ and $P_0$ are half-saturation constants. 

The Monod formulation guarantees strict resource limitation ($g'(N) > 0$), fulfilling assumption \eqref{hyp:H1}. Holling Type III response ensures that per-capita grazing pressure vanishes at low phytoplankton densities ($\lim_{P\to 0^+} h(t,P)/P = 0$), satisfying assumption \eqref{hyp:H2}. As the zooplankton metabolic rates $\mu_Z$ and $m_Z$ decouple from the threshold $\gamma_\text{crit}$, their parametrization is not required in our experimental setup.

%Phytoplankton persistence depends on the boundary dynamics where biomass approaches zero. Since the grazing functional response is of class $\mathcal{C}^1$ with $h(t,0)=0$, predator pressure vanishes at the origin. This decoupling isolates the phytoplankton persistence from the zooplankton parameters $\mu_Z$ and $m_Z$. Thus, defining the functional parameterizations for zooplankton grazing and mortality rates is unnecessary for the scope of this work. This independence represents a methodological advantage, as zooplankton grazing formulations typically introduce the highest parameter uncertainty in marine biogeochemical models \cite{Rohr2023}.

We adopted the thermodynamic framework proposed by Arroyo \etal \cite{arroyo} to model phytoplankton metabolic rates. This framework generate concave and convex Thermal Performance Curves (TPCs) that capture physiological non-linearities and thermal denaturation. The metabolic rate function, normalized to a reference temperature $T_0=\qty{293.15}{\kelvin}$ (\qty{20}{\degreeCelsius}), is expressed as
\begin{equation}\label{eqn:exp-eep}
k(T)=k(T_0)\exp\left[\frac{\Delta C}{R}\ln{\left(\frac{T}{T_0}\right)}+\frac{\Delta H}{R}\left(\frac{1}{T_0}-\frac{1}{T}\right)\right]\,,
\end{equation}
where $k(T_0)$ is the metabolic rate at the reference temperature, $R$ is the universal gas constant, $\Delta C$ is the change in heat capacity of activation, and $\Delta H$ is the enthalpy of activation. 

%Unlike the monotonic exponential formulations (e.g., $Q_{10}$ or Arrhenius) standard in classical NPZ models \cite{ROY201216, Gregg2008}, this formulation offers two distinct advantages: theoretical consistency, as it derives directly from EEP rather than empirical fits, and parsimony, requiring only two thermodynamic parameters ($\Delta C$ and $\Delta H$) which minimizes overfitting and computational cost in global simulations. To ensure ecological consistency, phytoplankton is modeled as a concave ectotherm subject to the thermodynamic constraints defined by Arroyo \etal \cite{arroyo}.

\subsection{Phytoplankton Parameters}

We parametrize a generic phytoplankton functional type with a maximum growth rate of $\mu_P=\qty{0.8}{\per\day}$ at \qty{20}{\degreeCelsius} \cite{Eppley1972, bissinger2008predicting}. The TPC follows a concave trajectory characterized by the following mean thermodynamic constraints defined by Arroyo \etal \cite{arroyo}
\begin{equation}
\Delta C_{P} = -14418.93~\mathrm{J\,mol^{-1}\,K^{-1}}, \qquad \Delta H_{P} = 4213943.43~\mathrm{J\,mol^{-1}}.
\end{equation}
Mortality is set to $m_P=\qty{0.024}{\per\day}$ at $20^\circ$C \cite{Baker2021}. The temperature-dependent rates are computed at each time step based on the local temperature profile. Table \ref{tab:model_parameters} summarizes the variables and constants used in the model.
%\begin{figure}[tb]
%    \centering
%    \includegraphics[width=0.9\textwidth]{figs/mu_p_metabolic_curve.pdf}
%    \caption{[Change figure]Thermal Performance Curve (TPC) for the maximum phytoplankton growth rate ($\mu_P$) derived from the thermodynamic model (\ref{eqn:exp-eep}). The solid teal line represents the statistical mean thermal response used for the global simulations. The background gray lines show an ensemble generated by randomly sampling $\Delta C$ and $\Delta H$ with a variance scaled to $0.025$ SD of the global empirical database \cite{arroyo}, illustrating the distribution of thermodynamic parameters around the mean baseline}
%    \label{fig:tpc_curve}
%\end{figure}

\begin{table}[tb]
    \centering
    \caption{Definition of model variables and parameters.}
    \label{tab:model_parameters}
    \begin{tabular}{lll}
        \toprule
        \textbf{Symbol} & \textbf{Description} & \textbf{Units} \\
        \midrule
        $H(t)$          & Mixed layer depth                    & \unit{\metre} \\
        $T(t)$          & Sea Surface Temperature              & $\mathrm{K}$ \\
        $K_d$           & Light attenuation coefficient        & \unit{\per\metre} \\
        $I(t)$          & Surface irradiance                   & \unit{\micro\mole\per\metre\squared\per\second} \\
        $N_0, P_0$      & Half-saturation constants            & \unit{\mmol\per\metre\cubed} \\
        $\mu_P$         & Max. phytoplankton growth rate       & \unit{\per\day} \\
        $m_P$           & Phytoplankton mortality rate         & \unit{\per\day} \\
        \bottomrule
    \end{tabular}
\end{table}

\section{Results: Global Baseline of Phytoplankton Annual Cycles}
\subsection{The Critical Nutrient Baseline}
\begin{figure}[tb]
    \centering
    \includegraphics[width=\textwidth]{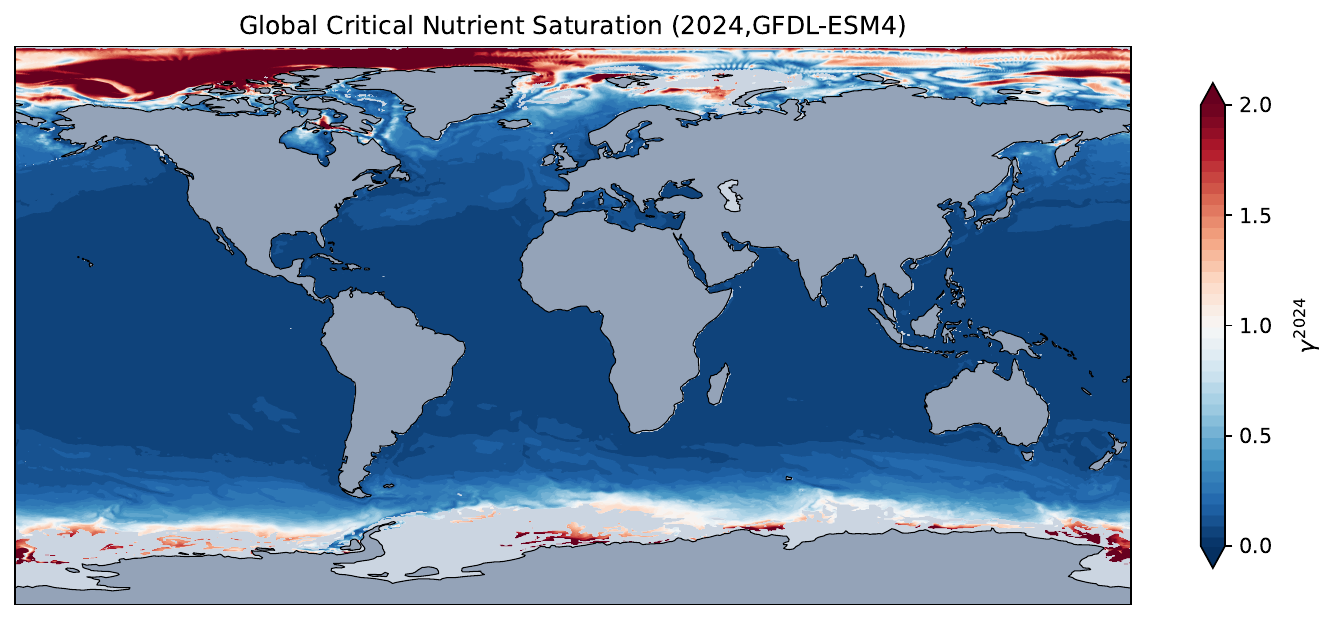}
    \caption{ Global Map of the Critical Nutrient Saturation($\gamma_\text{crit}$) derived from GFDL-ESM4 model data (2024).}
    \label{fig:gamma_crit}
\end{figure}
Figure~\ref{fig:gamma_crit} shows the global ocean distribution of present (2024) critical nutrient saturation. It reveals a clear spatial gradient. In the subtropics, low $\gamma_{\text{crit}}$ values indicate a low-cost regime, where phytoplankton populations can recover and sustain their annual cycles even under oligotrophic conditions. Conversely, high-latitude regions operate near the theoretical limit ($\gamma_{\text{crit}} = 1$), meaning they require near-maximal nutrient saturation to offset high metabolic and physical losses.

\subsection{Geographical Dominance of Physical and Metabolic Drivers}

To compute the local partial derivatives of $\gamma_\text{crit}$ with respect to the environmental drivers, we generated a 3D perturbation grid (for temperature, mixing, and light amplitude) at each spatial cell. A multiple linear regression was then fitted to the resulting $\gamma_\text{crit}$ response surface to extract the sensitivity gradients. The $R^2$ metric measures the goodness-of-fit of this linear approximation against the exact non-linear threshold responses. For both contemporary (2024) and future (2100) climate scenarios, this first-order Taylor decomposition remains robust, yielding a global mean $R^2 =\GFDLHyperMeanRsq$ with over \GFDLHyperAboveThresholdPerc\% of the ocean area exceeding an $R^2 >\GFDLHyperThreshold$. 

To identify dominant ecosystem controls, we classify each region according to its primary driver
\begin{itemize}
    \item \emph{Thermal Metabolism Driven:} Regions where temperature-dependent kinetic dominate ($\mathcal{I}_T = \max\{\mathcal{I}_T, \mathcal{I}_H, \mathcal{I}_I\}$).
    \item \emph{Physical Mixing Driven:} Regions dominated by MLD-driven light limitation, dilution and resource vertical transport ($\mathcal{I}_H = \max\{\mathcal{I}_T, \mathcal{I}_H, \mathcal{I}_I\}$).
    \item \emph{Light Amplitude Driven:} Regions constrained primarily by shifts in seasonal solar radiation forcing ($\mathcal{I}_I = \max\{\mathcal{I}_T, \mathcal{I}_H, \mathcal{I}_I\}$).
\end{itemize}

\begin{figure}[tbp]
    \centering
    \includegraphics[width=\textwidth]{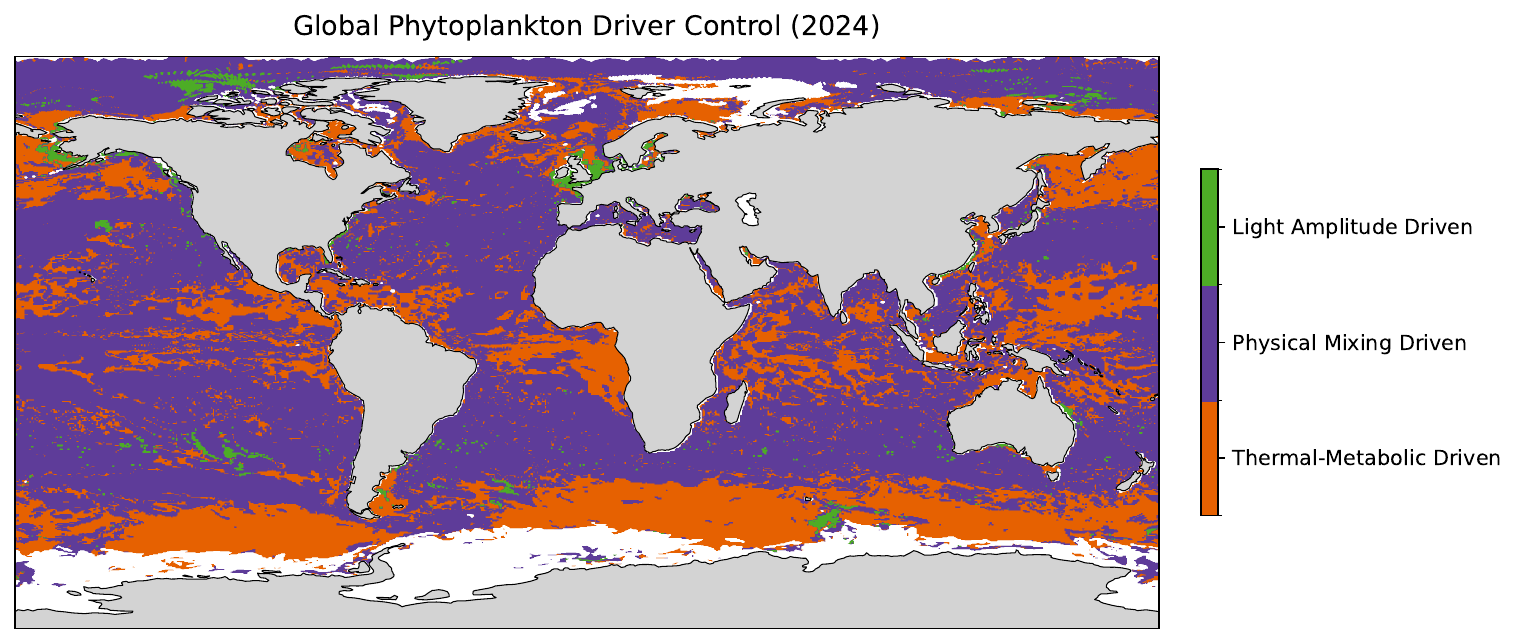}
    \caption{Climate drivers governing the phytoplankton annual cycle threshold ($\gamma_\text{crit}$) for the baseline year 2024. The map displays the dominant discrete regimes geographically (Orange: Thermal-Metabolic, Purple: Physical Mixing, Green: Light Amplitude). The contemporary ocean displays widespread dominance of physical mixing and thermal-metabolic drivers, while light amplitude control remains strictly localized.}
    \label{fig:drivers_baseline_2024}
\end{figure}

We note that according to Figure \ref{fig:drivers_baseline_2024}, there is no global primary driver in the 2024 baseline; however, the main drivers are physical-mixing and thermal-metabolic, while light-amplitude is confined.

\section{Global Projections: The Thermal Metabolic Shift}\label{sec:global_experiments}

We evaluate the long-term response of phytoplankton annual cycles under the SSP5-8.5 high-emission scenario.  Our choice is motivated by the objective of characterizing the potential impacts to seasonal recovery dynamics under extreme anthropogenic forcing.

\begin{figure}[tbp]
    \centering
    \begin{subfigure}[b]{\textwidth}
        \centering
        \includegraphics[width=\textwidth]{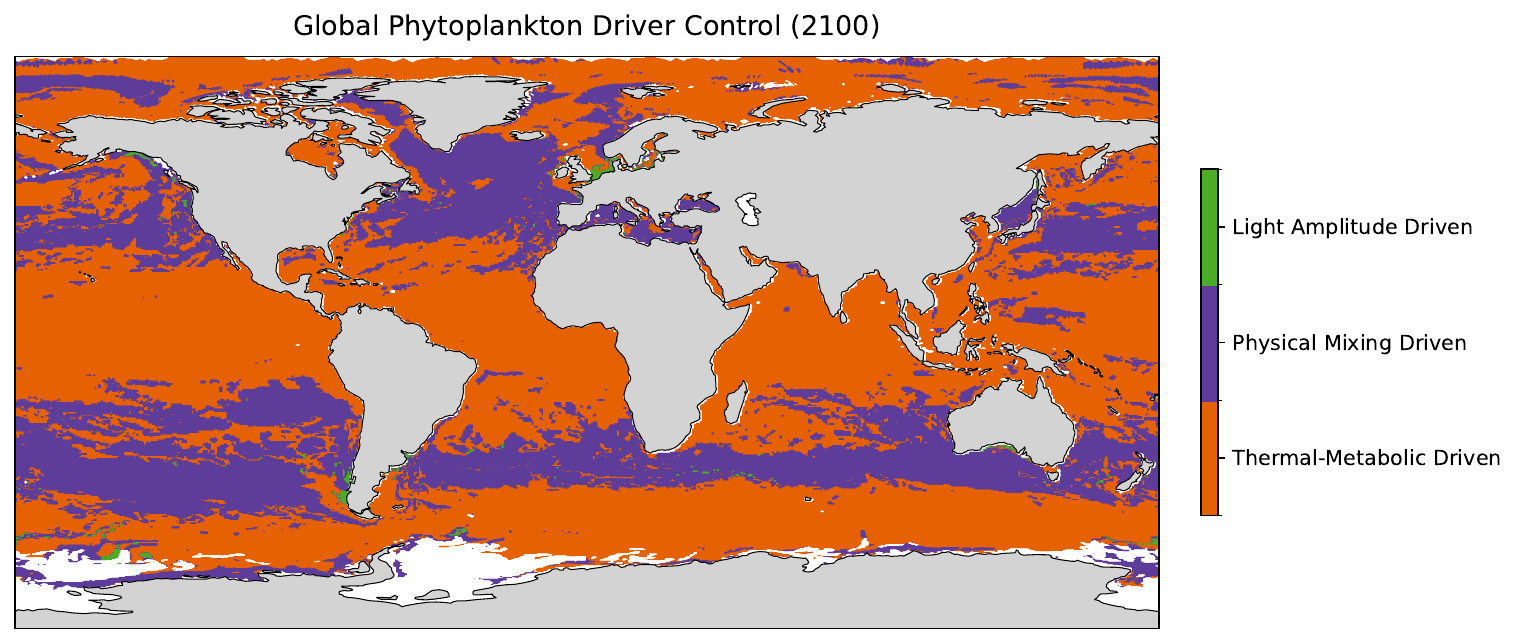}
        \caption{Discrete Regimes (2100)}
        \label{fig:discrete_2100}
    \end{subfigure}
    
    \vspace{0.5cm} % Espacio extra para que respiren los paneles
    
    \begin{subfigure}[b]{0.85\textwidth}
        \centering
        \includegraphics[width=\textwidth]{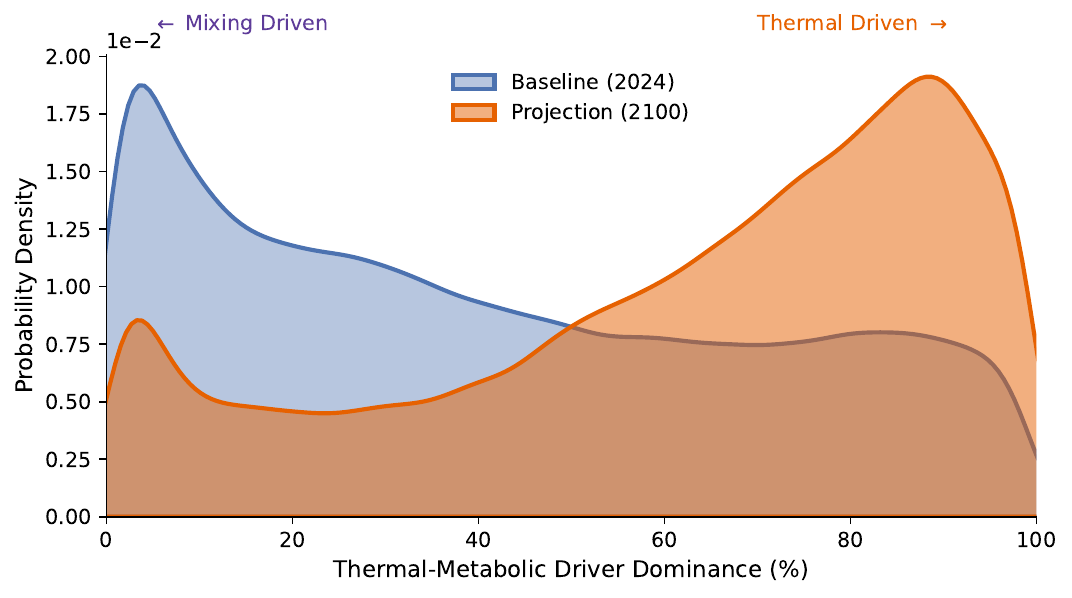}
        \caption{Global Driver Shift (2024 vs. 2100)}
        \label{fig:density_shift}
    \end{subfigure}
    \caption{Projected global distribution of climate drivers governing the phytoplankton annual cycle baseline ($\gamma_\text{crit}$) under the SSP5-8.5 scenario. Panel (a) maps the dominant discrete regimes geographically for 2100 (Orange: Thermal-Metabolic, Purple: Physical Mixing, Green: Light Amplitude). Panel (b) compares the 2024 and projected 2100 driver distributions using a 1D Kernel Density Estimate (KDE). The KDE illustrates the probability density of ocean pixels along the mixing-to-thermal dominance continuum, confirming a widespread metabolic tropicalization as the density peak shifts toward thermal dominance by the end of the century.}
    \label{fig:mechanisms_shift}
\end{figure}
According to Figure \ref{fig:mechanisms_shift} thermal metabolism emerges as the global primary driver of the 2100 projection. Compared to what we observe in the 2024 baseline, this projection is characterized by a homogenization where thermal metabolism displaces the historical influence of physical mixing across low and high latitudes. Besides this global thermalization trend, the North Atlantic Subpolar Gyre (SPG) stands out as a regional exception, maintaining its mixing-driven character.

%\subsection{The North Atlantic Subpolar Gyre (SPG): A Mixing-Driven Anomaly}

The SPG region is characterized by the ``warming hole'' phenomenon, persistent anomaly where SSTs remain nearly constant or even decrease deviating from the global warming trend \cite{rahmstorf_is_2024}. Under the high-emission forcing of the SSP5-8.5 scenario, as shown in Figure \ref{fig:mechanisms_shift}a, the mechanistic response of the SPG deviates from the global trend. While the global ocean shifts toward thermal-metabolic control by 2100, the SPG remains anchored in the physical mixing domain.

\section{Discussion}

\subsection{A General Approach for Sensitivity of Ecosystem Dynamics to Model Assumptions}

Here we demonstrate an approach to test the impacts of modifications to the parameterization of ecosystems in a climate model that requires minimal computational expense. By leveraging the linearization around the extinction equilibrium ($E_0$), we can analytically derive stability metrics such as the critical nutrient saturation($\gamma_\text{crit}$) and the trivariate driver decomposition. Earth System Models (ESMs) face a fundamental trade-off between biological realism and computational feasibility. Increasing ecosystem complexity not only demands greater processing power but also introduces structural uncertainties that obscure the causal link between physical forcing and biological rates \cite{anderson_2005}, a challenge that persists even in the latest CMIP6 generation where increased complexity has not necessarily reduced projection uncertainty \cite{kwiatkowski, heneghan_2021}.

This approach provides a reproducible tool for rapidly evaluating impacts of the climate change on ecosystems in a way that incorporates novel understanding of ecosystem processes \cite{martin_when_2024}. Whereas standard outputs (e.g., chlorophyll) typically show the net result of competing stressors, our metrics explicitly disentangle the mechanistic contributions of thermal-metabolic, physical-mixing and light-amplitude constraints. This attribution is important for identifying regional anomalies, such as the North Atlantic ``warming hole'', where the physical drivers diverge from global trends.

\subsection{Floquet Theory for Ecosystem Persistence}

By applying Poincaré maps, we compress continuous intra-annual environmental fluctuations into a discrete dynamical system. We establish uniform persistence for a periodically forced plankton system based on general functional properties, proving that weak occasional recovery mathematically guarantees a permanent annual cycle population (Proof of Theorem \ref{th:persistence}).

Theoretical ecology utilizes Floquet theory to evaluate species persistence and invasion fitness in seasonal environments \cite{klausmeier_floquet_2008}. Evaluating mutant invasion fitness typically requires intensive numerical routines to calculate Floquet multipliers over the resident seasonal cycle \cite{best_2023, ferris_2018}. Our framework bridges the gap for system persistence, while resolving coexistence equilibria for invasion fitness presents further challenges. By decoupling the Jacobian matrix at the extinction equilibrium via hypothesis \eqref{hyp:H2}, we evaluate the local Floquet exponent $\lambda_P$ algebraically. Thus, we demonstrate that this local linear instability translates to uniform persistence.

This provides theoretical guarantees for the persistence of annual cycles, without the need for numerical integration. Furthermore, it establishes a mathematical baseline to extend this framework toward eco-evolutionary invasion criteria.

\subsection{Limitations and Future Directions}

While our framework isolates the primary physical and metabolic controls on phytoplankton annual cycles, the assumptions inherent in our modeling choices suggest avenues for future research.

First, we distinguish between biological traits and biological rates. In this study, the physiological traits defining the TPC, such as the optimal temperatures, were held constant to represent a generic plankton functional type. Thus, our results isolate the thermodynamic pressure exerted by the environment from the compensatory mechanism of evolutionary adaptation. Future work should relax the assumption of fixed traits by incorporating adaptive dynamics \cite{Litchman2007,walworth_microbial_2020}. This would allow us to evaluate whether trait adaptation can keep pace with environmental warming, or if evolutionary lags will ultimately disrupt the annual cycle and drive local extinctions.

Second, our current analysis focuses on the annual cycles of primary producers. However, zooplankton grazing parameters, such as handling time and search efficiency, are also temperature-dependent \cite{ZHAO2020}. Integrating these thermally sensitive functional responses is a next step to evaluate the annual cycles of higher trophic levels and their long-term coexistence. Expanding our framework to fully coupled seasonal dynamics would allow us to identify regions where rapid warming triggers a metabolic mismatch, potentially breaking the seasonal synchronization between producers and consumers. 

Finally, to ensure a global attractor, the state space was defined as a compact set. While this closure is a necessary idealization to formally isolate the internal seasonal cycle from external noise, oceanic ecosystems are open systems. Future theoretical extensions could relax this compactness to analyze how non-autonomous open-boundary fluxes might destabilize the ecological attractor.

\section{Conclusion}

This study establishes a non-autonomous framework extending the Critical Depth Hypothesis (CDH) to evaluate the physical and metabolic drivers governing phytoplankton annual cycles. By analyzing invasion dynamics at the boundary, we derive the critical nutrient saturation($\gamma_\text{crit}$) as a structural baseline required to sustain seasonal population recovery under environmental forcing.

Applied to 21st-century climate projections from the GFDL-ESM4 model under the SSP5-8.5 scenario, our framework reveals a global regime shift toward thermal-metabolic dominance. This shift is widespread across most latitudes, with the notable exception of the North Atlantic Subpolar Gyre (SPG). Protected by anomalous thermal stasis, the SPG stands out as a mixing-driven anomaly where deep seasonal stratification continues to dictate the annual cycle, effectively resisting the global thermalization trend.

Beyond mapping these geographic transitions, our trivariate decomposition provides a computationally efficient approach to disentangle competing thermal, physical, and light forcings directly from ESM physics. By directly linking physical forcing with biological constraints, this framework offers a robust approach to project how climate-driven shifts in ocean mixing and thermal regimes will reshape the baseline of marine ecosystems in the coming century.

\appendix

\section{Mathematical Foundations: Floquet Theory and Uniform Persistence}\label{sec:foundations}

We analyze the stability of the NPZ-t model, as the coefficients (light, temperature and mixed layer depth) satisfy $A(t)=A(t+T_p)$ with $T_p=365~[\unit{days}]$,  autonomous stability is insufficient.

\subsection{Local Stability and the Monodromy Matrix}

The stability of the extinction equilibrium $E_0 = (C_0,0,0)$ -- where both Phytoplankton ($P$) and Zooplankton ($Z$) are absent -- is determined by the linearization of the system around the origin. We define the monodromy matrix $C$ as the operator that maps the initial perturbation $\Phi(0)=(P(0),Z(0))$ to its state after one full cycle
\begin{equation*}
\Phi(T_p) = C \Phi(0)\,.
\end{equation*}

The eigenvalues of $C$, known as Floquet multipliers ($\rho$), allow us to define the spectral radius $r(C) = \max |\rho|$, that dictates the local fate of the system
\begin{itemize}
\item If $r(C) < 1$, $E_0$ is locally asymptotically stable, \emph{i.\,e.}, any small population perturbation in a neighborhood $U \subset \mathbb{R}^2_+$ centered at $E_0$ will decay to zero ($\lim_{t \to \infty} \Phi(t) = 0$).
\item If $r(C) > 1$, $E_0$ is unstable.
\end{itemize}

\subsection{From Instability to Uniform Persistence}

While local instability suggests growth near the origin, the long-term survival of the phytoplankton population requires uniform persistence. This means there exists an $\epsilon > 0$ such that for any initial condition with $P(0)> 0$, $Z(0)\geq 0$
\begin{equation*}
\liminf_{t \to \infty} P(t) \geq \epsilon
\end{equation*}

Our analytical framework (Theorem \ref{th:persistence}) leverages the fact that for this periodic NPZ structure, the instability of $E_0$ (characterized by a positive Floquet exponent $\lambda_P = \frac{1}{T_p} \ln r(C) > 0$) is sufficient to guarantee uniform persistence for the phytoplankton population. In this regime, the extinction state acts as a repeller, ensuring the ecosystem remains bounded away from the boundary of the state space.

Ecologically, this mathematical boundary condition guarantees that the phytoplankton population avoids extinction and successfully closes its annual cycle, establishing the baseline threshold discussed in the main text.

\subsection{Key References}

\begin{itemize}
\item \emph{Floquet Theory:} \cite{Hale} for linear periodic systems and monodromy matrices and \cite{klausmeier_floquet_2008} for useful examples applied to ecology.
\item \emph{Persistence in Periodic Models:} \cite{zhao, wang2008threshold} for the formal link between the spectral radius of the periodic operator and global ecosystem persistence.
\end{itemize}

\section{Proofs of Theorems}\label{sec:proofs}

\subsection{Model Well-posedness}

\counterwithin{equation}{section}

\begin{proposition}\label{th:well_pos}
The system \eqref{eqn:NPZ-t-full} admits a unique global solution. The non-negative orthant $\mathbb{R}^3_+$ is positively invariant. Moreover, the total mass $C(t) = N(t) + P(t) + Z(t)$ is conserved for all $t \geq 0$, and the solutions are confined to the compact simplex 
\begin{equation*}
X = \{ (N,P,Z) \in \mathbb{R}^3_+ : N + P + Z = C_0 \}\,.
\end{equation*}
\end{proposition}

To establish the well-posedness of the model, we verify the existence, uniqueness, positivity, and conservation of mass of the solutions. Let $\Phi = (N, P, Z)^\top \in \mathbb{R}^3$ be the state vector. The system \eqref{eqn:NPZ-t-full} can be written as the Cauchy problem
\begin{equation} \label{eqn:cauchy}
    \frac{d \Phi}{d t} = F(t, \Phi), \quad \Phi(0) = \Phi_0 \in \mathbb{R}^3_+\,.
\end{equation}
The functions involved in the vector field $F(t, \Phi)$ are $C^1$ and bounded on the domain of interest. Therefore, $F$ is locally Lipschitz continuous, and by the Picard-Lindelöf Theorem, there exists a unique local solution $\Phi(t)$ on $[0, T_{\max})$.

Regarding positivity, the planes $P=0$ and $Z=0$ are invariant sub-spaces of the system, since $P=0 \implies \dot{P}=0$ and $Z=0 \implies \dot{Z}=0$ (noting that $h(0)=0$). Also, on the boundary $N=0$, we have:
\begin{equation*}
\left. \frac{d N}{d t} \right|_{N=0} = (1-\alpha) \mu_Z(t)h(P) Z + m_P(t) P + m_Z(t) Z^2 + s(t)(P+Z) \ge 0\,,
\end{equation*}
for all $(P,Z) \in \mathbb{R}^2_+$. Thus, the non-negative orthant $\mathbb{R}^3_+$ is positively invariant.

Adding the equations, we find $\frac{d}{dt}(N+P+Z) = 0$, which implies $C(t) = C_0$ for all $t \ge 0$. This conservation law, combined with the positivity of the variables, ensures that the solutions are confined to the compact simplex:
\begin{equation*}
X = \{ (N,P,Z) \in \mathbb{R}^3_+ : N+P+Z = C_0 \}\,.
\end{equation*}
Hence, the solutions are bounded, and finite-time blow-up is impossible, proving that $T_{\max} = \infty$. \qed

\subsection{Proof of Theorem \ref{th:extinction}}

We establish global extinction conditions using standard comparison arguments. First, we construct a differential inequality for the phytoplankton equation $P(t)$. From assumption \eqref{hyp:H1}, we know that the nutrient limitation term satisfies $g(N) \le g(C_0)$ forall $N\in [0,C_0]$. So, we obtain
\begin{align*}
    \frac{d P}{dt} & = f(t)\mu_P(t)g(N) P - \mu_Z(t)h(P) Z - m_P(t) P - s(t) P \\
                     & \le \left( f(t)\mu_P(t)g(C_0) - m_P(t) - s(t) \right) P\,.
\end{align*}
Let $P_l(t)$ be the solution to the following linear comparison system
\begin{equation*}
    \frac{d P_l}{dt} = \left( f(t)\mu_P(t)g(C_0) - m_P(t) - s(t) \right) P_l\,.
\end{equation*}
A standard comparison argument shows that as $P(0) = P_l(0)$, then $0 \le P(t) \le P_l(t)$ for all $t\ge0$.
Since the coefficients are $T_p\text{-periodic}$ we can write 
\begin{equation*}
    P_l(nT_p)=P(0)e^{nT_p \lambda_P}.
\end{equation*}
As $n\to \infty$, $P_l(nT_p)\to 0$, since $ \lambda_P < 0$. The per-capita growth rate of $P_l$ is bounded from above, letting it be $M$
\begin{equation*}
\frac{\dot{P_l}}{P_l}\leq M\,.
\end{equation*}
Let $n\in \mathbb{N}$, the above implies that for all $t\in\left[nT_p,nT_p+T_p\right]$ we get
\begin{equation*}
P_l(t)\leq P_l(nT_p)\exp(M(t-nT_p))\leq P_l(nT_p)\exp(MT_p)\,.
\end{equation*}
As $n\to \infty$, implying that $\lim_{t \to \infty} P_l(t) = 0$. Thus
\begin{equation*}
    \lim_{t \to \infty} P(t) = 0\,.
\end{equation*}

Recall the dynamics of $Z$ 
\begin{equation*}
    \frac{d Z}{dt} = \alpha \mu_Z(t) h(P) Z - m_Z(t) Z^2 - s(t) Z\,.
\end{equation*}
Since $h(P)$ is continuous, $\mu_Z$ is bounded and $h(0)=0$, the fact that $P(t)\to 0$ as $t\to \infty$ implies that for any $\eta > 0$, there exists a time $t_0$ such that $\alpha \mu_Z(t) h(P(t)) < \eta$ for all $t > t_0$. Choosing $\eta < \tfrac{1}{T_p}\int_{0}^{T_p}s(\tau)d\tau$, we obtain the inequality
\begin{equation*}
    \frac{d Z}{dt} \le \left( \eta - s(t) \right) Z -m_Z(t) Z ^2 \le \left( \eta - s(t) \right) Z\,,
\end{equation*}
Integrating this inequality, we have
\begin{equation*}
Z(t) \le Z(t_0) \exp \left( \int_{t_0}^t (\eta - s(\tau)) d\tau \right)=Z(t_0)F(t),\quad t\ge t_0\,.
\end{equation*}
As $F(t_0+T_p)<1$, due to the choice of $\eta$, and by a similar argument as in the case of phytoplankton, this implies that 
\begin{equation*}
    \lim_{t \to \infty} Z(t) = 0\,.
\end{equation*}

Since $P(t) \to 0$ and $Z(t) \to 0$ as $t\to \infty$, by conservation of the mass (Proposition \ref{th:well_pos}), the nutrient concentration must approach the total inventory $N(t) \to C_0$. Thus, the system globally converges to the extinction equilibrium $E_0 = (C_0, 0, 0)$.
\qed

\subsection{Proof of Theorem \ref{th:persistence}}

We will use Poincaré map theory applied to non-autonomous periodic dynamical systems. First recall the state space of the system $X:=\{(N,P,Z)\in \mathbb{R}^3:N+P+Z=C_0\}$ given by Proposition \ref{th:well_pos}, that is a compact simplex. Then we define the Poincaré map $f:X\to X$ given by $f(x_0)=\Phi(T_p,x_0)$, where $\Phi$ comes from the Cauchy problem \eqref{eqn:cauchy}. The map $f$ is continuous on $X$ due to the continuous dependence of the solutions with respect to the initial conditions. Also, we define $X_0:=\{(N,P,Z)\in X: P>0\}$ as the phytoplankton persistence space.

The main result is equivalent to that the Poincaré map of the $P$ component is uniformly persistent with respect to $(X_0, \partial X_0)$, \emph{i.\,e.}, 
\begin{equation*}
    \exists \eta>0, \quad \forall x \in X_0, \quad \liminf_{n\to\infty} P_{x}(nT_p) \geq \eta.
\end{equation*}
First, we note that $f(X_0)\subset X_0$. This is obtained directly from 
\begin{equation*}
P_{n+1}=P_n \exp\left( \int_0^{T_p} \mathcal{H}(\Phi(t,x_n))dt\right):=P_n G(x_n)\,,
\end{equation*}
where $x_n=f^n(x_0)=(N_n,P_n,Z_n)$ is the state of the system at the $n-$th period, and
\begin{equation*}
\mathcal{H}(t,N,P,Z)=f(t)\mu_P(t)g(N)-\mu_Z(t)\frac{h(P)}{P}Z-m_P(t)-s(t)
\end{equation*}
is the per-capita growth rate of the phytoplankton. The continuity of $\mathcal{H}$ near $E_0$ follows directly from assumption \eqref{hyp:H2}, which guarantees that the limit $\lim_{P\to 0^+} h(P)/P = h_0$ is well-defined and finite. By the analytical expression of the solution $P$, any trajectory with $P(0)>0$ remains strictly positive for any finite time, so $f^n(X_0)\subset X_0 \quad \forall n\in \mathbb{N}$. Then, given that the state space $X$ is compact and forward-invariant under $f$, it follows that $f$ is compact and point dissipative. Thus, by Theorem 1.1.3 from \cite{zhao}, there exists a global attractor $A\subset X$. 

Let us study the behavior near the extinction equilibrium $E_0$. Since $\lambda_P>0$, it follows that $G(E_0)=\exp(T_p\lambda_P)>1$. Let $\sigma$ be a constant such that $1<\sigma<G(E_0)$. By the continuity of $\mathcal{H}$ and the flow, there exists a neighborhood $\mathbb{B}_{\delta}(E_0)$ with $\delta>0$ such that $G(x)>\sigma$ for all $x\in \mathbb{B}_{\delta}(E_0)$. 

The map $f$ is weakly uniformly persistent with respect to $(X_0, \partial X_0)$, if and only if, 
\begin{equation*}
    \exists \eta_w>0, \quad \forall x \in X_0, \quad \limsup_{n\to\infty} P_{x}(nT_p) \geq \eta_w
\end{equation*}
Suppose that $f$ is not weakly uniformly persistent. Thus, for any $\eta_w > 0$, there exists an initial condition $x \in X_0$ such that
\begin{equation*}
    \limsup_{n\to\infty} P_x(nT_p) < \eta_w
\end{equation*}
By choosing $0 < \eta_w < \delta/3$ sufficiently small, there exists an initial condition $x \in X_0$ such that $P_n < \eta_w$ for all $n \ge N_1$. Since $\dot{P}/P \le M$, this yields $P(t) < \eta_w \exp(M T_p)$ for all $t \ge N_1 T_p$, which by the comparison argument in Theorem \ref{th:extinction} forces $Z_x(t) \to 0$ as $t \to \infty$. By mass conservation ($N_n = C_0 - P_n - Z_n$), $x_n$ eventually remains strictly inside $\mathbb{B}_{\delta}(E_0)$. For some index $N_0$, the sequence $(P_n)_{n \ge N_0}$ remains strictly inside $\mathbb{B}_{\delta}(E_0)$, yielding
\[
P_{N_0+k}=P_{N_0+k-1}G(x_{N_0+k-1})\geq P_{N_0+k-1}\sigma\geq \dots \geq P_{N_0}\sigma^k\,.
\]
As $\sigma>1$ and $P_{N_0}>0$, the sequence $P_{N_0+k}\geq P_{N_0}\sigma^k$ diverges to infinity as $k\to \infty$. This contradicts the fact that $X$ is compact invariant set (specifically $P(t)\leq C_0 \quad \forall t\geq 0$). Thus, $f$ is weakly uniformly persistent with respect to $(X_0,\partial X_0)$. Theorem 1.3.3 from \cite{zhao} implies that $f$ is uniformly persistent, \emph{i.\,e.}, there exists $\eta_u>0$ such that
\begin{equation*}
\forall x\in X_0, \quad \liminf_{n\to \infty}{P_{x}(nT_p)}\geq \eta_u.
\end{equation*}

Finally, we know that the maximum per-capita death rate is finite. Letting it be $L$ we obtain
\begin{equation*}
\frac{\dot{P}}{P}\geq -L\,.
\end{equation*}
Let $n\in \mathbb{N}$, the above implies that for all $t\in\left[nT_p,nT_p+T_p\right]$ we get
\begin{equation*}
P(t)\geq P(nT_p)\exp(-L(t-nT_p))\geq P(nT_p)\exp(-LT_p)\geq \eta_u \exp(-LT_p)=\eta >0\,.
\end{equation*}
The above holds true every year, therefore, we see that persistence is also uniform in the continuous system, \emph{i.\,e.},
\begin{equation*}
\liminf_{t\to \infty}{P_{x}(t)}\geq \eta \quad \forall x\in X_0\,.
\end{equation*}
\qed

%Following a similar argument to Theorem \ref{th:extinction}, since the per-capita growth rate of the phytoplankton is bounded from above, the convergence of the discrete sequence $P_{x_0}(nT_p)\to 0$ implies that $P_{x_0}(t)\to 0$ as $t\to \infty$. Thus, the zooplankton population also collapses ($Z_{x_0}(t)\to 0$), meaning that the $\omega$-limit set of any such trajectory is the singleton $\{E_0\}$.

%Consequently, as the trajectory approaches $E_0$ (where $Z \to 0$), the per-capita grazing loss term $\mu_Z(t)\frac{h(P)}{P}Z$ vanishes continuously.

%By convergence to $E_0$, there exists some $N_0$ such that the trajectory $(P_{N_0+n})_n$ remains inside the ball $\mathbb{B}_{\delta}(E_0)$ for all $n \ge 0$, therefore
%\[
%P_{N_0+k}=P_{N_0+k-1}G(x_{N_0+k-1})\geq P_{N_0+k-1}\sigma\geq \dots \geq P_{N_0}\sigma^k\,.
%\]
%As $\sigma>1$ and $P_{N_0}>0$, the sequence $P_{N_0+k}\geq P_{N_0}\sigma^k$ diverges to infinity as $k\to \infty$. This contradicts the fact that $X$ is compact invariant set (specifically $P(t)\leq C_0 \quad \forall t\geq 0$) and also contradicts that $x_n\to E_0$. Thus, $f$ is weakly uniformly persistent with respect to $(X_0,\partial X_0)$. Since $f$ has a global attractor and is weakly uniform persistent, Theorem 1.3.3 from \cite{zhao} implies that $f$ is uniformly persistent, \emph{i.\,e.}, there exists $\epsilon>0$ such that

\section{North Atlantic SPG Density Distribution}

\begin{figure}[tbp]
    \centering
    \includegraphics[width=0.85\textwidth]{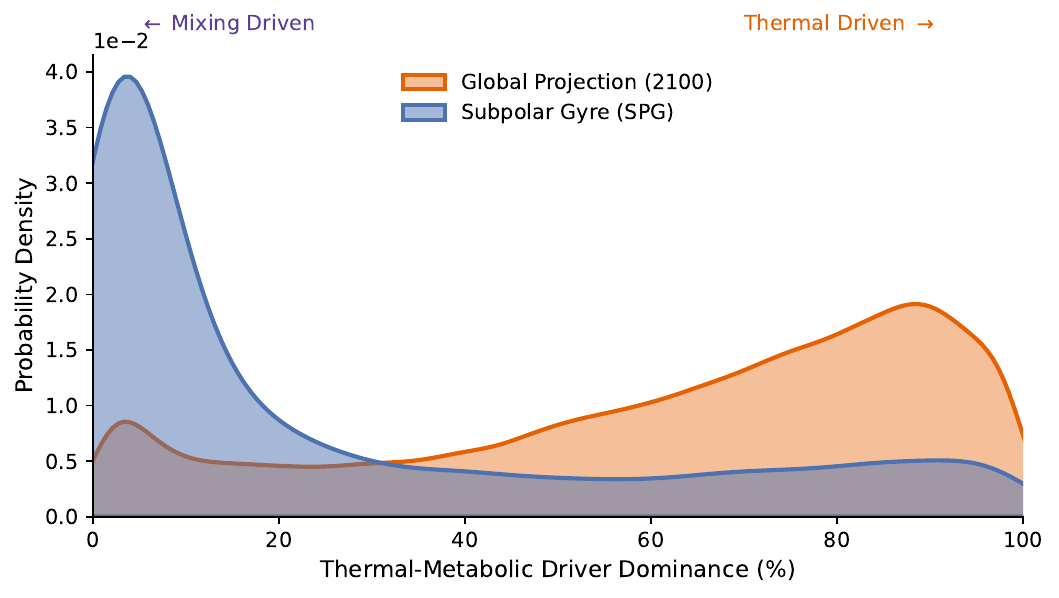} 
    \caption{Phytoplankton driver dominance in the North Atlantic Subpolar Gyre (SPG) compared to the global ocean for the year 2100 under the SSP5-8.5 scenario. The 1D Kernel Density Estimate (KDE) illustrates the probability density of ocean pixels along the mixing-to-thermal dominance continuum. While the global ocean is projected to shift toward thermal-metabolic control, the density distribution of the SPG remains anchored in the physical mixing domain. This contrast highlights the SPG as a regional exception, where seasonal stratification and dilution stress continue to govern the annual cycle baseline despite the global thermalization trend.}
    \label{fig:spg_analysis_2100}
\end{figure}

As shown in Figure \ref{fig:spg_analysis_2100}, the mechanistic response of the SPG deviates from the global trend. While the global ocean probability density shifts toward thermal-metabolic control by 2100, the SPG density distribution remains anchored in the physical mixing domain. 

\section*{Acknowledgments}

Work funded by the Franco-Chilean Binational Center of Artificial Intelligence, ANID Strengthening R\&D capabilities Program CTI230007 Inria Chile, and Inria Challenge OcéanIA (desc. num 14500). Pablo A. Marquet acknowledges support from Centro de Modelamiento Matem\'atico (CMM), Grant FB210005, BASAL funds for centers of excellence from ANID-Chile, and Proyecto Exploraci\'on 13220168. The authors are also grateful to Aurora Llanos and Ignacio Díaz for their valuable support during the final stages of this project; their assistance with technical adjustments and their careful proofreading were instrumental in the preparation of this manuscript. During the preparation of this work, the authors used Generative
AI technologies solely to edit and polish the English language and improve sentence structure.

\section*{Author Contributions}

All authors have made substantial intellectual contributions to the conception, execution, or design of the work. All authors have read and approved the final manuscript. Contributions, as specified by the Contributor Role Taxonomy (CRediT) can be summarized as:
\begin{itemize}
    \item Conceptualization: Matías Neto, Pablo A. Marquet, Mara Freilich, Luis Martí, Nayat Sanchez-Pi.
    \item Data curation: Mara Freilich, Matías Neto, Luis Martí.
    \item Formal Analysis: Matias Neto.
    \item Funding acquisition: Nayat Sanchez-Pi.
    \item Investigation:  Matías Neto, Pablo A. Marquet, Mara Freilich, Luis Martí, Nayat Sanchez-Pi.
    \item Methodology: Matías Neto, Pablo A. Marquet, Mara Freilich.
    \item Project administration: Nayat Sanchez-Pi.
    \item Resources: Matías Neto, Pablo A. Marquet, Mara Freilich, Luis Martí, Nayat Sanchez-Pi.
    \item Software: Matías Neto, Luis Martí.
    \item Supervision: Pablo A. Marquet, Nayat Sanchez-Pi.
    \item Validation: Matías Neto, Pablo A. Marquet, Mara Freilich, Luis Martí, Nayat Sanchez-Pi.
    \item Visualization: Matías Neto, Luis Martí.
    \item Writing—Original Draft Preparation: Matías Neto.
    \item Writing—Critical Review and Editing: Pablo A. Marquet, Mara Freilich, Luis Martí, Nayat Sanchez-Pi.
\end{itemize}
% All authors agree to be responsible for all aspects of the work, ensuring that any questions related to accuracy or integrity are investigated and resolved appropriately. A large language model (LLM) based on artificial intelligence was used to edit and polish the written text in order to improve grammar and overall style.
% \subsection{Conflicts of Interest}
% The authors declare there are no conflicts of interest
% \subsection{Ethics Statement}
% The authors confirm that this study is based on theoretical mathematical modeling and numerical simulations using publicly available climate model data. Therefore, the work did not involve human participants, animals, or primary data collection that would require review by an institutional review board or ethics committee.

\section*{Access to code} 
Code is available at \url{https://github.com/Inria-Chile/NPZ-t_DriverControl}.

\bibliographystyle{siamplain}
\bibliography{npz-t}

\end{document}